\documentclass[12pt]{amsart}
\usepackage{amsmath}
\usepackage{amsthm}
\usepackage{amssymb}

\usepackage{eucal}
\usepackage{times}
\usepackage{euler}

\usepackage{hyperref}

\author {Beno\^\i{}t Collins}
\address{Research Institute for Mathematical Sciences,
Kyoto university, Kyoto 606-8502, Japan}
\email{collins@math.kyoto-u.ac.jp}

\author{Piotr \'Sniady}
\address{Institute of Mathematics,
University of Wroclaw, \newline pl.\ Grunwaldzki~2/4, 50-384
Wroclaw, Poland} \email{Piotr.Sniady@math.uni.wroc.pl}

\title{New scaling of Itzykson-Zuber integrals}

\theoremstyle{plain}
\newtheorem{lemma}{Lemma}
\newtheorem{theorem}[lemma]{Theorem}

\newtheorem{hypothesis}[lemma]{Hypothesis}

\theoremstyle{definition}
\newtheorem*{definition}{Definition}

\theoremstyle{remark}

\newcommand{\R}{{\mathbb{R}}}

\newcommand{\gwia}{^{\star}}

\newcommand{\M}[1]{M_{#1}(\mathbb{C})}

 \DeclareMathOperator{\Tr}{Tr}

\DeclareMathOperator{\diag}{diag}

\begin{document}

\begin{abstract}
We study asymptotics of the Itzykson-Zuber integrals in the scaling
when one of the matrices has a small rank compared to the full rank.
We show that the result is basically the same as in the case when
one of the matrices has a fixed rank. In this way we extend the
recent results of Guionnet and Ma\"\i{}da who showed that for a
latter scaling the Itzykson-Zuber integral is given in terms of the
Voiculescu's $R$--transform of the full rank matrix.
\end{abstract}

\maketitle

\section{Introduction and the main result}

For diagonal matrices $A_N,B_N\in\M{N}$ we consider the
Iztykson--Zuber integral
$$ I^{(\beta)}_N(A_N,B_N)=\int e^{N
\Tr U A_N U\gwia B_N}  dm_N^{\beta}(U),$$ where $m_N^{\beta}$
denotes the Haar measure on the orthogonal group $O_N$ when
$\beta=1$, on the unitary group $U_N$ when $\beta=2$, and on the
symplectic group $Sp(N/2)$ when $\beta =4$ (in the latter case, $N$
is even). Usually one is interested in the study of the asymptotics
of the Itzykson-Zuber integral as the size $N$ of the matrices tends
to infinity and for this reason one would like to have an insight
into the limit
\begin{equation} \lim_{N\to\infty} \frac{1}{N^2} \log
I^{(\beta)}_N(A_N,B_N)=: \tilde{I}^{(\beta)}(\mu_A,\mu_B)
\label{eq:IZ}
\end{equation}
when the spectral measures $\mu_{A_N}, \mu_{B_N}$ converge weakly
towards some probability measures on $\R$. Existence of the limit of
\eqref{eq:IZ} was proved by Guionnet and Zeitouni
\cite{GuionnetZeitouni2002spherical}.

In \cite{Collins2002}, the study of another scaling of
Itzykson-Zuber integrals was initiated, namely when the rank $M$ of
the matrix $A_N$ is small compared to the full rank $N$. In order to
obtain a non-trivial limit one should consider rather the limit
\begin{equation} \lim_{N\to\infty} \frac{1}{MN} \log
I^{(\beta)}_N(A_N,B_N). \label{eq:IZ2}
\end{equation}

In this paper we shall say that a family of matrices $B_{N}\in\M{N}$
\emph{converges in distribution} if and only if for all integer $k$,
$N^{-1}\Tr (B_{N}^{k})$ tends towards a finite limit as
$N\to\infty$.

The first author proved that if $A_N=\diag(t,0,0,\dots)$ has rank
one and spectral measures of $B_N$ converge to some probability
measure $\mu_B$ then
$$\lim_{N\to\infty} \frac{1}{N} \frac{\partial^k}{\partial t^k} \log
I^{(2)}_N(A_N,B_N)\Big|_{t=0} = \frac{d^k}{dt^k}\int_0^t
R_{\mu_{B}}(s) ds \Big|_{t=0}$$ The function $R$ is called
Voiculescu's $R$-transform and is of central importance in free
probability theory. We recall its definition in Section
\ref{preuve}.

In this paper, the following hypothesis and notation about the
sequence of matrices $A_{N}\in\M{N}$ will be used frequently.
\begin{hypothesis}\label{hypo1}
$A_{N}$ is diagonal and it has rank $M(N)=o(N)$. Its non-zero
eigenvalues are denoted by $a_{1,N}\geq \ldots \geq a_{M,N}$.
\end{hypothesis}

Using the results presented in our previous article
\cite{CollinsSniady2004} it is possible to prove that the following
quite general statement holds true:
\begin{theorem}
Assume that $B_{N}$ has a limiting distribution $\mu_{B}$
and that $A_{N}$ is
uniformly bounded and satisfies Hypothesis \ref{hypo1}. Then for
$\beta=1$ and $2$,
\begin{multline*}
\Big| \frac{\partial^k}{\partial t^k}
\frac{1}{NM(N)}\log I_{N}^{(\beta)}(tA_{N},B_{N})\\
-\frac{\beta}{2M(N)}\sum_{i=1}^{M(N)}\frac{d^k}{dt^k}
\int_0^{\frac{2t}{\beta}} R_{\mu}(ta_{i,N}) ds \Big|_{t=0}=o(1).
\end{multline*}
If $\beta=4$ then
\begin{multline*}
\Big| \frac{\partial^k}{\partial t^k}
\frac{1}{NM(N)}\log I_{N}^{(\beta)}(tA_{N},B_{N})\\
+\frac{1}{M(N)}\sum_{i=1}^{M(N)}\frac{d^k}{dt^k} \int_0^{t}
R_{\mu}(-ta_{i,N}) ds \Big|_{t=0}=o(1).
\end{multline*}
\end{theorem}
Observe that it is convenient to separate the case $\beta=4$ from
the cases $\beta=1,2$.

We shall not prove this theorem since it is an immediate consequence
of Theorem 5.5 of \cite{CollinsSniady2004}. On the other hand, it
has been proved recently in Theorem 7 of \cite{GuionnetMaidaRtransform},
 that if the $L^{\infty}$-norm of $A_{N}$ is bounded by a constant depending on $B_{N}$
 (see Theorem \ref{GM} in this paper), under Hypothesis \ref{hypo-fond} for $B_{N}$,
and provided that
$M(N)=O(N^{1/2-\varepsilon})$ for some $\varepsilon >0$, one has for
$\beta=1,2$:
\begin{multline}\label{GMhigher}
\Big|\frac{1}{NM(N)}\log I_{N}^{(\beta)}(A_{N},B_{N})
-\frac{\beta}{2M(N)}\sum_{i=1}^{M(N)} \int_0^{\frac{2t}{\beta}}
R_{\mu}(s a_{i,N}) ds \Big|=o(1).
\end{multline}

The main result of this paper is Theorem \ref{res1}.
It shows that Equation \eqref{GMhigher}
still holds true if one replaces
restriction $M(N)=O(N^{1/2-\varepsilon})$ by Hypothesis \ref{hypo1},
and that it is enough to assume that the
$L^{\infty}$-norm of $A_{N}$ is uniformly bounded
by any constant, independently on $B_{N}$.

The techniques used in this paper are elementary, and substantially
simplify the proofs of
\cite{GuionnetMaidaRtransform}.
Better, it becomes possible to express the
$R$-transform as a limit Guionnet-Zeitouni's limits obtained in
\cite{GuionnetZeitouni2002spherical}.

We make the following hypothesis on the family of matrices
$(B_{N})_{N\in\mathbb{N}}$ with  $B_{N}\in\M{N}$:
\begin{hypothesis}\label{hypo-fond}
Let $b_{1,N}\geq\ldots \geq b_{N,N}$ be the real eigenvalues of
the %Hermitian
diagonal matrix
 $B_{N}$ and assume that the measures
 $N^{-1}\sum_{i=1}^{N}\delta_{b_{i,N}}$
 converge weakly towards a compactly supported
 probability measure $\mu_{B}$ and that
 $b_{1,N}$ (resp. $b_{N,N}$) converge towards
 the upper bound $\lambda_{\rm{max}}$ (resp. the lower bound
 $\lambda_{\rm{min}}$) of the support of
 $\mu_{B}$.
\end{hypothesis}

The Harisch-Chandra formula in the case when $\beta=2$ and when the
eigenvalues have no multiplicity tells that
\begin{equation*}
I_{N}^{(2)}(A_{N},B_{N})= \frac{\det (e^{N\cdot a_{iN}
b_{jN}})_{i,j=1}^N}{\Delta(a_{iN})\Delta(b_{iN})}
\end{equation*}
where $\Delta$ denotes Van der Monde determinant. In particular it
implies that the value of $I_{N}^{(2 )}$ depends only on the
eigenvalues of $A_{N},B_{N}$ therefore the assumption that
$A_{N},B_{N}$ are diagonal can be removed in the case $\beta=2$.
However a priori it cannot be removed in the other cases (although
it is an open question to the authors' knowledge).

We recall some definitions that can be found in
\cite{GuionnetMaidaRtransform}. For a probability measure $\mu_{B}$,
call $H_{\mu_{B}}$ its \emph{Hilbert transform}
$$H_{\mu_{B}}(z)=\frac{1}{z-\lambda}d\mu_{B}(\lambda ).$$
It is invertible under composition in a complex neighbourhood of
infinity; and let us call $z^{-1}+R_{\mu_{B}}$ its inverse.
$R_{\mu_{B}}$ is called the \emph{$R$-transform} of Voiculescu.

Let $H_{\rm{min}}=\lim_{z\uparrow\lambda_{\rm{min}}}H_{\mu_{B}}(z)$
and
$H_{\rm{max}}=\lim_{z\downarrow\lambda_{\rm{max}}}H_{\mu_{B}}(z)$.

\begin{definition}
For $\beta\in\{1,2,4\}$, a probability measure $\mu_{B}$ and a real
parameter $t$ we define $f^{(\beta )}_{\mu_{B}}(t)$ as follows.

If $\beta=1,2$, the function is
$$f^{(\beta )}_{\mu_{B}}(t)=t v(t)-
\frac{\beta}{2}\int \log \big(1+\frac{2}{\beta}t
v(t)-\frac{2}{\beta}t \lambda\big) d\mu_{E} (\lambda)$$ with
\begin{equation*}
v(t)=
\begin{cases}
R_{\mu_{B}}(\frac{2}{\beta}t)    &\text{if}\,\, H_{min}\leq\frac{2t}{\beta}\leq H_{\rm{max}},\\
\lambda_{\rm{max}}-\frac{\beta}{2t}   & \text{if}\,\, \frac{2t}{\beta}> H_{\rm{max}},\\
\lambda_{\rm{min}}-\frac{\beta}{2t}    &\text{if}\,\,
\frac{2t}{\beta}< H_{\rm{min}}.
\end{cases}
\end{equation*}
Observe that in the particular case $t\in [H_{min},H_{max}]$,
$$f^{(\beta )}_{\mu_{B}}(t)=\int_0^{\frac{2t}{\beta}}R_{\mu}(st) ds.$$
In the case $\beta=4$ the function is defined by
$$f^{(4)}_{\mu_{B}}(t)=-f^{(2)}_{\mu_{B}}(-t).$$
\end{definition}

One can prove that this function is continuous with respect to $t$
and also with respect to the metric on the measures given by
\begin{equation*}
\begin{split}
d(\mu,\mu')=|\lambda_{max}-\lambda_{max}'|+|\lambda_{min}-\lambda_{min}'|+
\sup \Big\{ \int f (d\mu-d\mu'), \\
||f||_{\infty}=1 \,\,\text{and f of Lipschitz constant}\,\,\leq 1
\Big\} .
\end{split}
\end{equation*}

\begin{theorem}[\cite{GuionnetMaidaRtransform}, Theorem 1.2 and 1.6]\label{GM}
 Assume that the family of matrices $B_{N}$ satisfies Hypothesis
 \ref{hypo-fond}.
 For $\beta=1,2$, for any real number $t$,
$$\lim_{N\to\infty}N^{-1}\log I_{N}(B_{N},{\rm diag \rm} (t,0,\ldots ,0))=f_{\mu_{B}}^{(\beta)}(t).$$
 \end{theorem}

 Note that although Guionnet and Ma\"\i{}da in \cite{GuionnetMaidaRtransform}
 do not handle the
case $\beta =4$, it is easy to extend their results to this case and
we will take these results for granted. We shall prove:
\begin{theorem}\label{res1}
Let $\beta\in\{1,2,4\}$.
\begin{itemize}
\item
Let $A_{N}$ be a sequence of matrices satisfying Hypothesis
\ref{hypo1}. Assume that $B_{N}$ satisfies Hypothesis
\ref{hypo-fond}. Assume furthermore that $\sup ||A_{N}||<\infty$
and $\beta=1,2,4$. Then
\begin{equation*}
\Big| \frac{1}{NM(N)}\log I_{N}^{(\beta)}(A_{N},B_{N})
-\frac{\beta}{2M(N)}\sum_{i=1}^{M(N)} f_{\mu_{B}}^{(\beta)}(a_{i,N})
\Big|=o(1).
\end{equation*}
\item
Let $\mu,\nu$ be compactly supported real probability measures.
Assume that $\mu$ has a connected support
$[\lambda_{{\rm min \rm}},\lambda_{{\rm max \rm}}]$,
and that there exists a constant $c>0$ such that
$$c\cdot dx1_{[\lambda_{{\rm min \rm}},\lambda_{{\rm max \rm}}]}\leq \mu .$$
%can we remove the above inequality ?
For
$a\in [0,1]$, let $\nu_{a}=a\nu + (1-a)\delta_{0}$. Then
$$\lim_{a\to 0}a^{-1}\tilde{I}^{(\beta )}(\nu_{a},\mu)=
\int_{t\in\mathbb{R}}f_{\mu}^{(\beta)}(t)d\nu (t).$$
\end{itemize}
\end{theorem}

The next section is devoted to proving this theorem.

\section{Proof of the main result}\label{preuve}

$I_{N}^{(\beta)}$ satisfies the following obvious translation
invariance property:
\begin{equation}\label{invariance}
I_{N}^{(\beta)}(A,B)=e^{N^{2}x}I_{N}^{(\beta)}(A+x\rm{Id},B),
\end{equation}
therefore there is no loss of generality in assuming that $B_{N}$ is
positive.

We make temporarily the following assumption
\begin{hypothesis}\label{temporaire}
There exists a constant $c>0$ such that
$$b_{i+1,N}+c/N\geq b_{i,N}.$$
\end{hypothesis}

This implies that the distribution of $\mu_{B}$ has absolutely
continuous part with respect to Lebesgue measure whose density is
$\geq c$ almost everywhere on its support, and the support should be
an interval.

We shall need the following technical result:

\begin{lemma}\label{minimax}
Let $B\in\M{N}$ be Hermitian, with eigenvalues $b_{1}\geq\ldots \geq
b_{d}\geq 0$. Let $\Pi$ be a projector of rank $N-1$ and $B'=\Pi
B\Pi$ and let $b_{1}'\geq\ldots \geq b_{d}'=0$ its eigenvalues. Then
one has
$$b_{1}\geq b_{1}'\geq b_{2}\geq\ldots\geq b_{d}\geq b_{d}'=0.$$
\end{lemma}

\begin{proof}
Let $V$ be the set of Hermitian projections of rank $n+1-i$ and $V'$
be the subset of projections of $V$ dominating $1-\Pi$. According to
the ``minimax'' theorem, $b_{i}=\min_{\pi\in V}||\pi B \pi||$, and
$b_{i-1}'=\min_{\pi\in V'}||\pi B \pi||$. Since $V'\subset V$ this
implies that $b_{i+1}'\geq b_{i}$. Replacing $B$ by $||B||Id-B$
shows that $b_{i}\geq b_{i}'$, which concludes the proof.
\end{proof}

\begin{lemma}\label{positivite}
Let $A,A',B\in\M{N}$. If $A\geq A' $ and $B\geq 0$ then
$$I_{N}^{(\beta)}(A,B)\geq I_{N}^{(\beta)}(A',B).$$
\end{lemma}

\begin{proof}
Observe that if $A,B\geq 0$ then $\Tr AB\geq 0$. Therefore $\Tr
AUBU^{*}\geq \Tr A'UBU^{*}\geq 0$, which implies $\exp (N\Tr
AUBU^{*})\geq \exp (N\Tr A'UBU^{*})$. Integrating over the Haar
measure on the unitary group yields the desired result.
\end{proof}

\begin{proof}[Proof of Theorem \ref{res1}]

In the case $\beta=2$, we identify $U_{N-1}$ with a subgroup of
$U_{N}$ via the morphism
$$U_{N-1}\to
\left(\begin{array}{cc}1 & 0 \\0 & U_{N-1}\end{array}\right)\subset
U_{N}.$$

In the case $\beta=1$ we consider instead the embedding of
$O_{N-1}$ into $O_{N}$
$$O_{N-1}\to
\left(\begin{array}{cc}1 & 0 \\0 & O_{N-1}\end{array}\right)\subset
O_{N},$$ and in the case $\beta=4$, the embedding $SP_{N-1}$ into
$SP_{N}$
$$SP_{N-1}\to
\left(\begin{array}{cc}Id_{2} & 0 \\0 &
SP_{N-1}\end{array}\right)\subset SP_{N}.$$

We shall only deal the case $\beta=2$, the cases $\beta=1,4$ being
similar.

Through the above identification, and using the invariance of Haar
measure by convolution, we have
\begin{multline}\label{11}
I^{(\beta)}_N(A,B)=\\
\int_{U\in U_{N}}\int_{V\in U_{N-1}} \exp (N\Tr
AVUBU^{*}V^{*})dm_{N}^{\beta}(U)dm_{N-1}^{\beta}(V).
\end{multline}

In the remainder of the proof, $N$ is fixed, so in order to lighten
the notation we
omit the subscript $N$ for the matrices $A_{N}$ and $B_{N}$.
Since a unitary conjugation of $A$ by a unitary leaves
$I_{N}^{(\beta)}(A,B)$ invariant, one can assume without loss of
generality that
$$A={\rm diag\rm} (a_{1},\ldots ,a_{N})=
\left(\begin{array}{cc}a_1 & 0 \\0 & A_{-}\end{array}\right)=
\left(\begin{array}{cc}A_{+} & 0 \\0 & a_{N}\end{array}\right)$$
with $a_{i}\geq a_{j}$ if $i<j$. We shall also adopt the following
notation:
$$\tilde{B}=U BU^{*} =
\left(\begin{array}{cc}\tilde{b}_{11} & \tilde{B}_{1-} \\
\tilde{B}_{-1} & \tilde{B}_{--}\end{array}\right).
$$
Observe that $\tilde{B}_{--}\in\M{N-1}$ and and
$\tilde{b}_{11}\in\mathbb{C}$ are random matrix depending on the
matrix random variable $U$. With these notations we have
\begin{multline}\label{33}
I^{(\beta)}_N(A_N,B_N)=\\
\int_{U\in U_{N}}\int_{V\in U_{N-1}}
\exp (N\Tr A_{N}V\tilde{B}V^{*})dm_{N}^{\beta}(U)dm_{N-1}^{\beta}(V)\\
=\int_{U\in U_{N}}\int_{V\in U_{N-1}}
\exp (Na_{1}\tilde{b}_{11}+N\Tr A_{-}V\tilde{B}_{--}V^{*})dm_{N}^{\beta}(U)dm_{N-1}^{\beta}(V)\\
=\int_{U\in U_{N}} \exp (Na_{1}\tilde{b}_{11}) \int_{V\in U_{N-1}}
\exp (N\Tr A_{-}V\tilde{B}_{--}V^{*})dm_{N}^{\beta}(U)dm_{N-1}^{\beta}(V)\\
=\int_{U\in U_{N}} \exp (Na_{1}\tilde{b}_{11})
I^{(\beta)}_{N-1}(\frac{N}{N-1}A_{-},\tilde{B}_{--})dm_{N}^{\beta}(U).
\end{multline}

Let $B^{+,i}={\rm diag \rm}(b_{1},\ldots ,b_{N-i})$ and $B^{-,i}=
{\rm diag \rm}(b_{i},\ldots ,b_{N})$. By Lemma \ref{minimax} and Lemma
\ref{positivite}, one has
\begin{multline}\label{22}
I^{(\beta)}_{N-1}\Big(\frac{N}{N-1}A_{-},B^{-,1}\Big)\leq\int_{U\in
U_{N}}
I^{(\beta)}_{N-1}\Big(\frac{N}{N-1}A_{-},\tilde{B}_{--}\Big)\\
\leq I^{(\beta)}_{N-1}\Big(\frac{N}{N-1}A_{-},B^{+,1}\Big).
\end{multline}
Inequality \eqref{22} applied to the last line of Equation
\eqref{33} implies that
\begin{multline}\label{44}
I^{(\beta)}_{N-1}\Big(\frac{N}{N-1}A_{-},B^{-,1}\Big) \int_{U\in
U_{N}}
\exp (Na_{1}\tilde{b}_{11})dm_{N}^{\beta}(U)\\
\leq \int_{U\in U_{N}} \exp (Na_{1}\tilde{b}_{11})
I^{(\beta)}_{N-1}\Big(\frac{N}{N-1}A_{-},\tilde{B}_{--}\Big)dm_{N}^{\beta}(U)\\
\leq
 I^{(\beta)}_{N-1}\Big(\frac{N}{N-1}A_{-},B^{+,1}\Big)
\int_{U\in U_{N}} \exp (Na_{1}\tilde{b}_{11})dm_{N}^{\beta}(U)
\end{multline}
In addition, it is obvious that
$$\int_{U\in U_{N}}
\exp (Na_{1}\tilde{b}_{11})dm_{N}^{\beta}(U)
=I^{(\beta)}_N({\rm diag \rm} (a_{1},0,\ldots ,0),B_N).$$ Thus,
inequality \eqref{44}  and equation \eqref{33} imply
\begin{multline}\label{66}
I^{(\beta)}_{N-1}\Big(\frac{N}{N-1}A_{-},B^{-,1}\Big)
I^{(\beta)}_N({\rm diag \rm} (a_{1},0,\ldots ,0),B_N) \leq
I^{(\beta)}_N(A_N,B_N)\\
\leq I^{(\beta)}_N({\rm diag \rm} (a_{1},0,\ldots ,0),B_N)
I^{(\beta)}_{N-1}\Big(\frac{N}{N-1}A_{-},B^{+,1}\Big).
\end{multline}
A recursive application of bound \eqref{66} yields
\begin{multline}\label{55}
\prod_{i=1}^{M(N)}I_{N+1-i}\Big(\frac{N}{N+1-i}{\rm diag \rm} (a_{i},0,\ldots
,0), B^{-,i}\Big)
\leq I_{N}(A,B)\\
\leq
\prod_{i=1}^{M(N)}I_{N+1-i}\Big(\frac{N}{N+1-i}diag(a_{i},0,\ldots
,0), B^{+,i}\Big).
\end{multline}

Inequality \eqref{55} and change of formula or Equation
\eqref{invariance} and Theorem \ref{GM} imply Theorem \ref{res1}
under Hypothesis \ref{temporaire}.

Theorem \ref{res1} is proved in full generality without hypothesis
\ref{temporaire} by the continuity property of the limit
$f^{(\beta)}_{\mu}$ stated after its definition.

\end{proof}

\section{Concluding remarks}

In comparison to \cite{GuionnetMaidaRtransform}, our method uses
much weaker assumptions on the rank and norm of the matrices.
Better, it allows us to state a continuity result about the limit of
$IZ$ in a new scaling and thus validate an ``inversion of limit''
phenomenon.

Last, it shows that Hypothesis \ref{hypo-fond} is a relevant one to
perform computations outside the phase transition zone. This is a
substantial improvement to the paper \cite{GuionnetMaidaRtransform}
in which nothing was proved at ``high temperature'' for the finite
scaling of rank $M>1$.

Unfortunately, our approach heavily relies on real-number
inequalities and we believe that finer estimates, as well as
complex valued estimates of \cite{GuionnetMaidaRtransform} can not
be established with our methods.

\emph{Acknowledgements}: The authors benefited from an excellent
working atmosphere in Banff free probability meeting in October
2004. They also thank Alice Guionnet for stimulating discussions.

The research of P. \'S. was supported by State Committee for
Scientific Research (KBN) grant \mbox{2 P03A 007 23}, RTN network:
QP-Applications contract No.~HPRN-CT-2002-00279, and KBN-DAAD
project 36/2003/2004. This author is a holder of a scholarship of
European Post-Doctoral Institute for Mathematical Sciences.

The research of B.C. was partly supported by a JSPS postdoctoral
fellowship and by RIMS.

\bibliographystyle{alpha}
\bibliography{biblio}

\end{document}